\input amstex
\documentstyle{amsppt}
\input bull-ppt
\keyedby{bull347e/PAZ}

\define\Z{{\Bbb{Z}}}
\define\C{{\Bbb{C}}}
\define\Q{{\Bbb{Q}}}

\define\FF{{\cal{F}}}

\topmatter
\cvol{28}
\cvolyear{1993}
\cmonth{January}
\cyear{1993}
\cvolno{1}
\cpgs{95-98}

\title Representations and $K$-theory\\
of Discrete Groups\endtitle
\author Alejandro Adem\endauthor
\address Department of Mathematics, University of Wisconsin,
Madison, ~Wisconsin~53706\endaddress
\ml adem\@math.wisc.edu\endml
\thanks Research partially supported by an NSF 
grant\endthanks
\date May 20, 1992\enddate
\abstract Let $\Gamma$ be a discrete group of finite
virtual cohomological dimension with certain finiteness 
conditions of the
type satisfied by arithmetic groups.  We define a 
representation ring
for $\Gamma$, determined on its elements of finite order, 
which is of finite
type.  Then we determine the contribution of this ring to 
the topological
$K$-theory $K^*(B\Gamma )$, obtaining an exact formula for 
the
difference in terms of the cohomology of the centralizers 
of elements of
finite order in $\Gamma$.\endabstract
\endtopmatter

\document
\heading 0. Introduction\endheading

Let $\Gamma$ denote a discrete group of finite virtual 
cohomological
dimension.  Examples of this type of group include finite 
groups,
arithmetic groups, and mapping class groups, making them 
an important class of
objects in both topology and algebra.  In particular, 
understanding the
classifying space $B\Gamma$ for such groups is a central 
problem in
algebraic topology.  Unfortunately, the cohomology 
$H^*(B\Gamma ,\Z )$ is
a very intractable object; consequently, there are few 
available calculations
(e.g., see \cite{So}).  In sufficiently high dimensions 
the cohomology is known to
depend only on the lattice $\FF$ of finite subgroups in 
$\Gamma$ \cite{B, 
F},
but in general this yields a complicated spectral sequence 
involving
the cohomology of the normalizers $N(S)$, $S\in\FF$.

In this note we outline an approach to understanding the 
r\^ole of
representations in the topology of $B\Gamma$ as was done 
in the case of
finite groups by Atiyah \cite{At}.  We define a 
representation ring determined
on the elements of {\it finite order} in $\Gamma$, which 
for a large class
of groups (including arithmetic groups) is of finite rank. 
 Then we
indicate to what degree the topological $K$-theory 
$K^*(B\Gamma )$ is
determined by these representations.  In fact, we provide 
a precise
description of the discrepancy in terms of the {\it 
rational cohomology}
of the centralizers
of elements of finite order in $\Gamma$.  Complete details 
will appear
elsewhere.
\heading 1. A reduced representation ring for 
$\Gamma$\endheading

From now on we will assume that $\Gamma$ has a {\it 
finite} number of
distinct conjugacy classes of elements of finite order and 
that their
centralizers are homologically finite.  These hypotheses 
are known to
hold, in particular, for arithmetic groups.
\dfn{Definition 1.1}
Let $V,W$ be two finite-dimensional
$\C\Gamma$-modules.  We say that $V$ is $\FF$-isomorphic 
to $W$ if
$V\Big|_S\cong W\Big|_S$ for all $S\in\FF$.
\enddfn
\dfn{Definition 1.2}
$R_\FF (\Gamma)$ is 
the Grothendieck group on
$\FF$-isomorphism classes of finite-dimensional 
$\C\Gamma$-modules.
\enddfn

We can, of course, also describe $R_\FF (\Gamma )$ as a 
quotient of the usual
representation ring $R(\Gamma )$.  Let $n(\Gamma )$ denote 
the number
of distinct conjugacy classes of elements of finite order 
in $\Gamma$.
Using character theory arguments, we prove
\thm{Proposition 1.3}
$R_\FF (\Gamma )$ is a 
commutative, unitary
ring, which as an abelian group is free of rank $n(\Gamma 
)$.
In particular, $\Gamma$ is torsion-free if and only if 
$R_\FF (\Gamma )
\cong\Z$.
\ethm

Similarly, if $\FF (p)$ denotes the family of all finite 
$p$-subgroups of
$\Gamma$ and $n_p(\Gamma )$ the number of distinct 
conjugacy classes
of elements of order a power of $p$ ($p$ prime), then 
$R_{\FF (p)}(\Gamma )$
can be defined and will be of rank $n_p(\Gamma )$.  The 
following
examples illustrate that these rings are readily 
computable from
subgroup data, unlike the cohomology.
\ex{Example 1.4}
$\Gamma = \hbox{SL}_2(\Z )$, $n(\Gamma)=8$, and
$$
R_\FF (\hbox{SL}_2(\Z ))\cong \Z [w]\Biggm/ w^8+w^6-w^2-1=0.
$$
\endex
\ex{Example 1.5} 
$\Gamma = \hbox{SL}_3(\Z )$, $n_2(\Gamma)=5$, and
$$
R_{\FF (2)}(\hbox{SL}_3(\Z ))\cong
\Z [\alpha_1,\alpha_2,\beta_1,\beta_2]\Bigg/
\matrix\format\l\\
\alpha^2_1=\alpha^2_2=1,\;\alpha_1\beta_1=\beta_1,\\
\beta^2_1=2(1+\alpha_1),\;\beta_2^2=2(1+\alpha_2),\\
\alpha_2\beta_2=\beta_2,\;\alpha_1\alpha_2=\alpha_1+
\alpha_2-1,\\
\alpha_1\beta_2=2\alpha_1+\beta_2-2,\\
\beta_1\beta_2=2\beta_1+2\beta_2-4,\\
\alpha_2\beta_1=2\alpha_2+\beta_1-2.
\endmatrix
$$
\endex
\heading 2. Contribution to $K$-Theory\endheading

For the sake of clarity of exposition, we work at  a fixed 
prime $p$;
let $K^*_p(\;\; )$ denote $p$-adic $K$-theory and $\C_p$ 
the completion
of the algebraic closure of $\Q_p$.  We choose a fixed 
normal subgroup
$\Gamma'\subseteq\Gamma$ such that $\Gamma'$ is 
torsion-free,
so $G=\Gamma /\Gamma'$ is finite.  If $\gamma\in\Gamma$, 
let $C(\gamma )$
denote its centralizer; then it can be expressed as an 
extension
$$
1\rightarrow C(\gamma )\cap\Gamma'\rightarrow C(\gamma 
)\rightarrow H_\gamma
\rightarrow 1
$$
where $|H_\gamma |<\infty$.  Our main result is the 
following.
\thm{Theorem 2.1}
Let $\Gamma$ be a discrete group of finite
v.c.d.\ satisfying our finiteness assumptions.  Then there 
is an exact
sequence
$$
0\rightarrow I_p\rightarrow K^*_p(B\Gamma )\otimes \C_p
{\buildrel {\varphi_p}\over\longrightarrow} R_{\FF 
(p)}(\Gamma )\otimes
\C_p\rightarrow 0
$$
where
$\varphi_p$ is a surjection of rings,
and we have an additive decomposition
$$
I_p\cong\bigoplus_{(\gamma )}\tilde K^*_p(B(C(\gamma 
)\cap\Gamma'))^{H_\gamma}
\otimes \C_p
$$
where the sum is taken over conjugacy classes 
of elements of order a finite power
of $p$.
\ethm
\thm{Corollary 2.2}
$$
K^*_p(B\Gamma )\otimes \C_p\cong R_{\FF (p)}(\Gamma 
)\otimes\C_p
$$
if and only if $\tilde H^*(BC(\gamma ),\Q)\equiv 0$ for 
every element
$\gamma\in\Gamma$ of order a power of $p$.
\ethm

The corollary follows from the fact that $I_p$ is 
determined by the
{\it cohomology} of $BC(\gamma )$; it is, of course, {\it 
independent} of
the choice of the extension.
\ex{Example 2.3}
$\Gamma =G_1{\displaystyle{\mathop *}_H} G_2$ is
an amalgamated product of finite groups.  Then
$$\eqalign{
&K^0_p(B\Gamma )\otimes\C_p\cong R_{\FF (p)}(\Gamma 
)\otimes\C_p\, ,\cr
&K^1_p(B\Gamma )\otimes\C_p\cong (\C_p)^{v_p(\Gamma )}\cr}
$$
where 
$$v_p(\Gamma )=n_p(\Gamma )-n_p(G_1)-n_p(G_2)+n_p(H)$$ 
represents the
total sum of $\dim_{\Q}\,H^1(BC(\gamma ),\Q)$
as $\gamma\in\Gamma$ ranges 
over conjugacy 
classes of elements of order a
power of $p$.
\endex
\ex{Example 2.4}
$\Gamma =\hbox{ SL}_3(\Z )$ and
$$
K^*_2(B\Gamma )\otimes\C_2\cong R_{\FF (2)}(\Gamma 
)\otimes \C_2\, ,
$$
whence we can use Example 1.5 
to determine this ring (compare with \cite{So, 
TY}).
\endex
\ex{Example 2.5}
$\Gamma =\hbox{ GL}_{p-1}(\Z )$, where $p$ is odd prime.
If $\roman{Cl}(p)=$ class number of $p$, then $R_{\FF 
(p)}(\Gamma )$ can be
computed from the extension
$$
0\rightarrow R_{\FF (p)}(\Gamma 
)\rightarrow\left(\bigoplus_{\roman{Cl} (p)}R(
\Z /p)\right)^\Delta \rightarrow\Z^{t(p)-1}\rightarrow 0
$$
where $\Delta =$ Galois group
and $t(p)=$ number of $\Delta$-orbits in the set
of ideal classes.  Hence $rk_\Z R_{\FF (p)}(\Gamma )=1+
\roman{Cl} (p)$, and in this
case
$$
I_p\cong \tilde K^*_p(B\Gamma')^G\otimes\C_p\oplus \left( 
\bigoplus_{\roman{Cl} (p)}
\tilde K^*_p\left( (S^1)^{(p-3)/2}\right)\otimes 
\C_p\right)\, .
$$
\endex
\demo{Sketch of Proof of \rm2.1}
From a theorem of Serre \cite S for the
class of groups we consider
that there exists a finite-dimensional $\Gamma$-complex
$X$ with {\it finite} isotropy, contractible fixed point 
sets, and $X/\Gamma$
of finite type.  By essentially identifying bundles that 
agree on finite
subgroups, we construct a surjection of rings
$$
K^*_\Gamma (X)\twoheadrightarrow R_\FF (\Gamma )\, .
$$
Next we identify $K^*_\Gamma (X)\cong K^*_G(X/\Gamma')$ 
($\Gamma'$, $G$ as
before) and use an {\it additive} decomposition for 
$K^*_G(X/\Gamma')\otimes\Bbb{C}$
obtained previously by the author \cite{A} to estimate the 
kernel of this ring map
in terms of the centralizers of elements of finite order 
in $\Gamma$.

The final technical step is to complete this map, as by 
the Atiyah-Segal
completion theorem, $K^*_G(X/\Gamma')^{\wedge} \cong 
K^*(B\Gamma )$
(at $IG\subseteq R(G))$.  Doing this {\it locally} leads 
to the statement
in Theorem 2.1.\qed
\enddemo

\enddocument